\documentclass[12pt]{article}
\usepackage{amsmath,amssymb,amsthm,graphics,color}

\begin{document}

\begin{center}
{\bf Tiling Lattices with Sublattices, II}
\end{center}

\begin{center}
David Feldman, University of New Hampshire \\
James Propp, University of Massachusetts Lowell \\
Sinai Robins, Nanyang Technological University
\end{center}

\noindent
{\sc Disclaimer}: This document is a place-holder for a longer and more 
carefully written article that will include more than just a Mirsky-Newman-style
proof of the generalized Mirsky-Newman theorem of [FPR].  However, we are 
posting this rough draft here since there has been interest in the question of 
whether Fourier methods are an essential component of the proof in [FPR].

\bigskip

Call {\em Cartesian} those subgroups of ${\bf Z}^d$  that have the
form $\prod_i a_i{\bf Z}$.  Using Fourier techniques, Feldman, Propp
and Robins [FPR] proved their generalization of the Mirsky-Newman
theorem:

\bigskip\noindent
{\bf Theorem}  ${\bf Z}^d$ admits no decomposition as the disjoint
union of cosets of distinct Cartesian subgroups.

Here we present a generating function proof in the spirit of
Mirsky-Newman's original.

Where Mirsky-Newman analyzes poles of 1-variable generating functions,
we shall need an appropriate multi-variable analog. So, given a
rational function $f$ of $k$ complex variables, call a point
$p\in {\bf C}^k$ an {\em order $d$ pole} provided:\\
\hspace*{.3 in} 1) the denominator of $f$ vanishes at $d$;\\
\hspace*{.3 in} 2) for almost every line $L$ (= affine real subspace)
through $p$, $|f|$ has restriction to $L$ continuous away from $p$ on
an interval $I$ containing $p$, and such that for any smooth parameter
$t$ on $L$ which vanishes at $p$, $|f|$ grows on the order of
$t^{-d}$.

{\bf Proof}  Suppose that disjoint cosets $T_j$ of Cartesian subgroups
union to ${\bf Z}^d$. Then the disjoint sets $R_j=T_j\cap {\bf N}^d$
union  to ${\bf N}^d$, each $R_j$ has a $d$-variable generator
function
$$G_j:=\prod_i \frac{{x_i}^{m_{j,i}}}{1-{x_i}^{n_{j,i}}}$$
and
$$S:=\sum_j G_j =   \prod_i \frac{1}{1-{x_i}}\ .$$

Each $G_j$ has a $d$-order pole at $p_j:=(e^{2\pi \sqrt{-1}
/n_{j,i}})_i$ but $S$ does not (except in the trivial case where
$R_1=N^d$).\footnote{Note that while $p_j$ does not sit isolated,
merely as a pole, from the poles of $G_j$, it does sit isolated from
the $G_j$'s finitely many poles of order $d$. Indeed, $G_j$  has
$d$-order poles at all points of the form $(e^{2\pi \sqrt{-1}
k_i/n_{j,i}})_i$.  }

The rest follows [FPR]: pick $j$ with $\prod_i n_{j,i}$ maximized.
Some $d$-order pole of at least one $G_{j'}$, $j'\not=j$ must cancel
$G_j$'s pole at $p_j$.  By the choice of $j$, we must mean the pole of
$p_{j'}$, and then $n_{j,i}=n_{j',i}$ for all $i$, by the choice of
$j$. And that makes $T_j$ a translate of $T_j'$.

\bigskip\noindent [FPR] D. Feldman, J. Propp, S. Robins, Tiling Lattices
with Sublattices, I, Discrete and Combinatorial Geometry, to appear;
{\tt arXiv:0905.0441}.

\end{document}